\documentclass[12pt,draft]{amsart}
\usepackage{amsmath,amsthm,latexsym,amscd,amsbsy,amssymb}

 \relax

\chardef\bslash=`\\ 

\makeatletter
\def\verbatim{\interlinepenalty\@M \@verbatim
  \leftskip\@totalleftmargin\advance\leftskip2pc
  \frenchspacing\@vobeyspaces \@xverbatim}
\makeatother
\newtheorem{thm}{Theorem}[section]
\newtheorem{cor}[thm]{Corollary}
\newtheorem{lem}[thm]{Lemma}
\newtheorem{pro}[thm]{Proposition}
\newtheorem{rem}[thm]{Remark}
\newtheorem{q}[thm]{Question}
\newtheorem{defin}[thm]{Definition}
\newtheorem{ex}[thm]{Example}

\theoremstyle{definition}

\theoremstyle{remark}

\numberwithin{equation}{section}

\font\f=msbm10

\begin{document}


\title[$4$-manifolds]{On some dimensional properties of $4$-manifolds}
\author{Alex Chigogidze}
\address{Department of Mathematics and Statistics,
University of Saskatche\-wan,
McLean Hall, 106 Wiggins Road, Saskatoon, SK, S7N 5E6,
Canada}
\email{chigogid@math.usask.ca}
\thanks{Both authors acknowledge the support
of their respective science foundations: NSERC and INTAS;
and the second author recalls with gratitude the
hospitality offered by the Department of Mathematics $\&$ Statistics
at the University of Saskatchewan.}

\author{V.~V.~Fedorchuk}
\address{Chair of General Topology and Geometry, Faculty
of Mechanics
and Mathematics, Moscow State University, Moscow, 119899, Russia}
\email{vitaly@gtopol.math.msu.su}


\begin{abstract}{It is shown, under the assumption of
Jensen's principle $\lozenge$, that if for a complex $L$
with $\left[ L \right] \geq \left[ S^{4}\right]$
there exists a metrizable compactum whose extension
dimension is $L$, then there exists a differentiable,
countably compact, perfectly normal and hereditarily
separable $4$-manifold whose extension dimension is also
$\left[ L\right]$.}
\end{abstract}

\maketitle
\markboth{A.~Chigogidze and V.~V.~Fedorchuk}{$4$-manifolds}


\section{Introduction}\label{S:intro}
It was shown in \cite{vvf1} that under the assumption of
Jensen's principle $\lozenge$ there exists a differentiable
$n$-manifold $M_{m}^{n}$, $n \geq 4$, of any given
Lebesgue dimension $m$ where $m > n$. This manifold is
countably compact, perfectly normal and hereditarily separable.
Under the same set-theoretical assumption $\lozenge$ for any
countable ordinal number $\alpha > 4$ there exists
\cite{vvf} a $4$-manifold $M_{\alpha}^{4}$
with $\operatorname{Ind}M_{\alpha}^{4} = \alpha$. \cite{vvf}
also contains examples of: (a) weakly infinite-dimensional
$4$-manifolds without the large inductive dimension 
and (b) strongly infinite-dimensional $4$-manifolds.
Recently it was shown \cite{cv} that for a given
{\em countable} complex $L$,
with $\left[ L \right] \geq \left[ S^{4}\right]$
and which serves as the extension
dimension of a metrizable compactum, there exists a 
differentiable $4$-manifold $M = M^{4,L}$
with $\operatorname{e-dim}M = \left[ L \right]$. It
should be emphasized that it is still
unknown whether the extension dimension of a metrizable
compactum is realized by a countable complex.
Below
we construct a differentiable $4$-manifold with similar properties
for {\em any}, not necessarily countable, complex.


\section{Preliminaries}\label{S:pre}
We recall that a subset $U \subseteq X$ of a space $X$
is {\em functionally open} in $X$ if there is a continuous
map $\varphi \colon X \to \text{\f R}$ such that
$U = \varphi^{-1}(\text{\f R}-\{ 0\} )$. Also, we say that $X$
is at most $n$-dimensional (and write $\dim X \leq n$)
if every finite functionally open cover $\mathcal U$ of
$X$ has a finite functionally open refinement
$\mathcal V$ of order $\leq n+1$. The latter means that
$\cap{\mathcal V}^{\prime} = \emptyset$ for any family
${\mathcal V}^{\prime} \subseteq {\mathcal U}$ consisting
of at least $n+2$ elements. 

For normal spaces this definition is equivalent to the usual
definition of Lebesgue dimension. The next statement is well known
(see, \cite{pears}, for example).
\begin{pro}\label{P:2.1}
For every space $X$ we have $\dim X = \dim \beta X$.
\end{pro}

We assume that reader is
familiar with notions of {\em a CW-complex}, {\em a simplicial
complex with the metric
topology} and an {\em absolute neighborhood
retract} in the category $\mathcal M$ of metrizable
spaces ($ANR({\mathcal M})$-space) (see, for instance,
\cite{fri}). In what
follows, by a {\em simplicial complex} we mean
any simplicial complex with the metric topology. Let us note here that
all simplices are assumed to be closed which implies that every finite
simplicial complex is compact.
By an $ANR$ we mean an $ANR({\mathcal M})$-space.

\begin{thm}\cite[Theorem 3.3.10]{fri}\label{T:2.3}
Every simplicial complex is an $ANR$.
\end{thm}

The next statement, which is a corollary of
\cite[Theorem 5.2.1]{fri}, allows us to consider
only simplicial complexes.
\begin{thm}\label{T:2.4}
Every CW-complex has a homotopy type of a simplicial
complex.
\end{thm}

\begin{defin}\label{D:2.5}
{\em Following \cite{coho} we say that a space $Z$ is an} 
absolute extensor of a normal space $X$ {\em and write
$Z \in AE(X)$ if for each closed
subspace $Y$ of $X$ any map $f \colon Y \to Z$ has an extension
$\bar{f} \colon X \to Z$.}
\end{defin}

The next statement is an immediate corollary of the above
definition.

\begin{pro}\label{C:2.7}
If $Z \in AE(X)$, $X$ is a normal space and $Y$ is a
closed in $X$, then $Z \in AE(Y)$.
\end{pro}

\begin{defin}\label{D:2.9}
{\em Let $X$ and $Z$ be normal spaces. Recall that $Z$ is an}
absolute neighborhood extensor of a space $X$
{\em (notation: $Z \in ANE(X)$) if for
every closed subspace $Y \subseteq X$ any map
$f \colon Y \to Z$ has an extension $\bar{f} \colon U \to Z$
where $U$ is a neighborhood of $Y$ in $X$.}
\end{defin}

\begin{pro}\label{P:2.10}
Let $X$ be a normal countably compact space and let
$L$ be a simplicial complex. Then $L \in ANE(X)$.
\end{pro}
\begin{proof}
Let $f \colon Y \to L$ be a map of a closed subset $Y \subseteq X$.
Since $L$ is metrizable, $f(Y)$ is compact. Hence $f(Y)$
is contained in some finite subcomplex $K \subseteq L$.
But every finite complex is an $ANE$ for any normal space. Thus,
there is an extension $\bar{f} \colon U \to K$ of $f$ defined on an open
neighborhood $U$ of $Y$ in $X$.
\end{proof}

\begin{pro}\label{P:2.11}
The following conditions are equivalent for every
countably compact normal space $X$ and every simplicial complex $L$:
\begin{itemize}
\item[(1)]
$L \in AE(X)$.
\item[(2)]
$L \in AE(\beta X)$.
\end{itemize}
\end{pro}
\begin{proof}
$(1) \Longrightarrow (2)$. By Definition \ref{D:2.5}, we
need to check that for every closed set $Y \subseteq
\beta X$ any map $f \colon Y \to L$ has an
extension $\bar{f} \colon \beta X \to L$.
By Proposition \ref{P:2.10}, there is an
extension $f_{1} \colon U \to L$, where
$U$ is a neighborhood of $Y$ in $\beta X$.
Let $U_{1}$ be a smaller neighborhood of $Y$
in $\beta X$ such that
$\operatorname{cl}_{\beta X}U_{1} \subseteq U$.
Set $F = X\cap\operatorname{cl}_{\beta X}U_{1}$
and let $f_{2} = f_{1}|F$. By condition (1),
there is an extension $\bar{f}_{2} \colon X \to L$.
As in the proof of Proposition \ref{P:2.10},
$\bar{f}_{2}(X)$ is contained in some finite
complex $K \subseteq L$. But as was noted above
$K$ is compact. Hence the map $\bar{f}_{2}$
can be extended to a map $\bar{f} \colon
\beta X \to K \subseteq L$. It remains to
show that $\bar{f}|Y = f$. But $\bar{f}|F =
f_{1}$. Hence, since $F$ is dense in
$\operatorname{cl}_{\beta X}U_{1}$, we have
$\bar{f}|\operatorname{cl}_{\beta X}U_{1} =
f_{1}$. On the other hand, $f_{1}|Y = f$.

$(2) \Longrightarrow (1)$. Let $Y$ be a closed subset
of $X$ and let $f \colon Y \to L$ be a map. Set
$F = \operatorname{cl}_{\beta X}Y$. Since $Y$ is closed
in a normal space $X$, $F = \beta Y$. Then $f$ can
be extended to a map $f_{1} \colon F \to L$ because
$f(Y)$ lies in some finite complex $K \subseteq L$.
Now, by condition $(2)$, the map $f_{1} \colon F \to L$
can be extended to a map $\bar{f}_{1} \colon \beta
X \to L$. It only remains to note that the map
$\bar{f} = \bar{f}_{1}|X$ extends $f$. Proposition
\ref{P:2.11} is proved.
\end{proof}

\begin{pro}\label{P:2.12}
Let $X$ be a countably compact normal space,
$F$ be its closed subset and $U = X - F$. Suppose
$L \in AE(F)$ and $L \in AE(Y)$ for every closed in $X$
set $Y \subseteq U$. Then $L \in AE(X)$.
\end{pro}
\begin{proof}
By Definition \ref{D:2.5}, we need to verify that for
every closed set $A \subseteq X$ any map
$f \colon A \to L$ has an extension $\bar{f}
\colon X \to L$. Let $f_{0} = f|(A\cap F)$. Since
$L \in AE(F)$, the map $f_{0}$ can be extended to a
map $\bar{f}_{0} \colon F \to L$. Define the map
$f_{1} \colon A \cup F \to L$ by letting $f_{1}|A = f$
and $f_{1}|F = \bar{f}_{0}$. Clearly, $f_{1}$ is continuous.
By Proposition \ref{P:2.10}, the map $f_{1}$ has an
extension $\bar{f}_{1}\colon V \to L$, where $V$ is a
neighborhood of $A\cup F$ in $X$. Take a neighborhood
$V_{1}$ of $A\cup F$ such that
$\operatorname{cl}(V_{1}) \subseteq V$ and let
$Y = X-V_{1}$, $Y_{1} = \operatorname{Bd}(V_{1})$,
$g = \bar{f}_{1}|Y_{1}$. Then $Y$ is closed in $X$ and
$Y_{1}$ is closed in $Y$. By condition $L \in AE(Y)$.
Hence, the map $g$ has an extension $\bar{g} \colon
Y \to L$. Finally, define a map $\bar{f} \colon X \to L$ by
letting
\[ \bar{f}|Y = g \;\;\text{and}\;\; \bar{f}|
\operatorname{cl}(V_{1}) = \bar{f}_{1} .\]
Evidently, $\bar{f}$ is well defined and continuous. It
is also clear that $\bar{f}|A = f$. Proposition
\ref{P:2.12} is proved.
\end{proof}

Next we define a relation $\leq$ on the class of all
simplicial complexes. Following \cite{dra1} we say that $K \leq L$
if for every normal countably compact space $X$ the
condition $K \in AE(X)$ implies the condition $L \in AE(X)$.
The relation $\leq$ is reflexive and transitive
and, consequently, it is a relation of preorder.
This preorder
induces the following equivalence relation:
\[ K \sim L \Longleftrightarrow K \leq L\;\;\text{and}\;\;
L \leq K .\]
For a simplicial complex $L$ by $\left[ L\right]$ we denote
the class of all complexes which are equivalent to $L$.
These classes $\left[ L\right]$ are called {\em extension types}.

\begin{rem}\label{R:2.13}
{\em Relation $L \in AE(X)$, preorder $\leq$ and extension
types $\left[ L\right]$ can be defined for different classes
of spaces $X$. A.~N.~Dranishnikov \cite{dra2} defined
relation $L \in AE(X)$ for the class ${\mathcal MLC}$
of all metrizable locally compact spaces. In \cite{dra1} he
defined this relation for the class ${\mathcal C}$ of all
compact Hausdorff spaces. One can define relation $\leq_{\sigma}$
and associated concepts for arbitrary class $\sigma$
of topological spaces. Let ${\mathcal MC}$ be the class
of all compact
metrizable spaces and ${\mathcal CC}$ be the class of all
normal countably compact spaces.}
\end{rem}

\begin{pro}\label{P:2.14}
For any simplicial complexes $K$ and $L$ the following conditions
are equivalent:
\begin{itemize}
\item[(1)]
$K \leq _{{\mathcal MC}} L$;
\item[(2)]
$K \leq _{{\mathcal C}} L$;
\item[(3)]
$K \leq _{{\mathcal CC}} L$;
\item[(4)]
$K \leq _{{\mathcal MLC}} L$;
\end{itemize}
\end{pro}
The equivalence $(1) \Longleftrightarrow (2)$ was
proved in \cite[Theorem 11]{dra3}. For the equivalence
$(2) \Longleftrightarrow (3)$ consult with
Proposition \ref{P:2.11}. As for the equivalence
$(1) \Longleftrightarrow (4)$ it follows from Theorem
\ref{T:2.21} and the next trivial statement.

\begin{pro}\label{P:2.15}
Let $X_{\alpha}$, $\alpha \in A$ be a discrete family of
normal spaces. Then for any simplicial complex $L$
\[ L \in AE\left(\oplus\{ X_{\alpha} \colon \alpha
\in A\}\right) \Longleftrightarrow L \in AE(X_{\alpha})
\;\;\text{for each}\;\; \alpha \in A .\]
\end{pro}

If $\sigma$ is a class of topological spaces, then
by $\text{\f E}_{\sigma}$ we denote the class of all
extension types of
all simplicial complexes generated by the relation
$\leq_{\sigma}$. In view of Proposition \ref{P:2.14} we shall
use a simpler notation: $\text{\f E}$ and $\leq$.

\begin{defin}{\em (see \cite{coho}\cite{dra1})}\label{D:2.16}
{\em Let $X$ be a countably compact normal space. Its}
extension dimension {\em $\operatorname{e-dim}X$ is
defined as the smallest extension type
$\left[ L\right]$ of simplicial complexes,
satisfying condition $L \in AE(X)$.}
\end{defin}

\begin{pro}\cite[\S 3, Proposition 1]{dra1}\label{P:2.17}
For any compactum $X$ there exists unique extension
type $\left[ L\right]$ such that
$\operatorname{e-dim}X = \left[ L\right]$.
\end{pro}

\begin{pro}\cite[\S 3, Proposition 2]{dra1}\label{P:2.18}
The correspondence $\operatorname{e-dim}$ maps the class
$\mathcal C$ epimorphically onto the class $\text{\f E}$.
\end{pro}
Propositions \ref{P:2.11} and \ref{P:2.17} yield

\begin{pro}\label{P:2.19}
For any normal countably compact space $X$ there exists
unique extension type $\left[ L\right]$ such that
$\operatorname{e-dim}X = \operatorname{e-dim}\beta X =
\left[ L\right]$.
\end{pro}

Propositions \ref{P:2.11} and \ref{P:2.18} yield
\begin{pro}\label{P:2.20}
The correspondence $\operatorname{e-dim}$ maps the
class ${\mathcal CC}$ epimorphically onto the class
$\text{\f E}$.
\end{pro}

\begin{thm}\label{T:2.21}
Suppose that a normal countably compact space $X$ is
the union of its closed subsets $X_{i}$, $i \in \omega$.
If $\operatorname{e-dim}X_{i} \leq \left[ L\right]$
for each $i \in \omega$, then $\operatorname{e-dim}X
\leq \left[ L\right]$.
\end{thm}
The proof of the above statement repeats the proof
(see, for instance, \cite{alpas}) of
classical countable sum theorem for Lebesgue dimension
$\dim$ for normal spaces by means of extension of
mappings into $S^{n}$. The main feature of the sphere $S^{n}$
exploited in that proof is $S^{n} \in ANE(X)$.
The corresponding
property $L \in ANE(X)$ in our case is guaranteed
by Proposition \ref{P:2.10}.

For further references we formulate just mentioned
description of the
Lebesgue dimension as a separate statement. Obviously it
provides the main link between the theory of Lebesgue dimension
and the theory of extension dimension.

\begin{thm}\label{T:2.22}
For any normal space $X$,
\[ \dim X \leq n \Longleftrightarrow
\operatorname{e-dim}X \leq \left[ S^{n}\right] .\]
\end{thm}

\section{On a realization of dimensional
types by manifolds}\label{S:manifolds}

We recall one result from \cite{vvf} in a more convenient
for us form.
\begin{thm}\label{T:3.1}
For an arbitrary metrizable compactum $C$, assuming
$\lozenge$, there exists a differentiable, countably
compact, perfectly normal, hereditarily separable
$4$-manifold $M_{C}^{4}$ such that
$\beta M_{C}^{4} -M_{C}^{4}$ is a metrizable compactum
homeomorphic to the disjoint sum of $C$ and some open
subset $U$ of the $3$-dimensional sphere $S^{3}$.
\end{thm}
The manifold $M_{C}^{4}$ is a manifold of type $M_{\varphi}^{4}$
from \cite{vvf}, where $\varphi =
\varphi_{C} \colon B^{4} \to B_{\varphi}^{4}$ is a quotient mapping,
defined on the closed ball $B^{4}$, with the following properties.

Let the sphere $S^{3}$ be the boundary of $B^{4}$. There exists a
closed set $A \subseteq S^{3}$ such that
\begin{itemize}
\item[(i)]
$A = \varphi^{-1}\varphi (A)$;
\item[(ii)]
$\varphi (A) = C$;
\item[(iii)]
each fiber $\varphi^{-1}(y)$, $y \in C$, is a non-degenerate
continuum nowhere dense in $S^{3}$.
\item[(iv)]
$\varphi^{-1}\varphi (x) = \{ x\}$ for every $x \in B^{4}-A$.
\end{itemize}

Thus, $\varphi (S^{3}) \equiv S_{\varphi}^{3}$ is homeomorphic to the
disjoint sum of $C$ and $S^{3}-A$. By \cite[Proposition 2.3]{vvf},
$\beta M_{C}^{4} - M_{C}^{4} = S_{\varphi}^{3}$.

Let $\Lambda$ be the class of all complexes and let
\[ \Lambda^{0} = \{ L \in \Lambda \colon \left[ L \right] =
\operatorname{e-dim} X\;\text{for some metrizable compactum} X\} .\]
By Proposition \ref{P:2.17}, for every metrizable
compactum $X$ there is a complex $L \in \Lambda^{0}$
such that $\operatorname{e-dim} X = \left[ L\right]$. Set
\[ \Lambda_{4}^{0} = \{ L \in \Lambda^{0} \colon
\left[ L\right] \geq \left[ S^{4}\right] \} .\]

The next theorem is the main result of this section.
\begin{thm}\label{T:3.2}
For an arbitrary complex $L \in \Lambda_{4}^{0}$, assuming
$\lozenge$, there exists a differentiable, countably compact,
perfectly normal, hereditarily separable $4$-manifold $M = M^{4,L}$
such that $\operatorname{e-dim} M = \left[ L\right]$.
\end{thm}
\begin{proof}
We use the scheme of the proof of \cite[Theorem 3.1]{cv},
where a similar result was obtained for countable complexes. The
only difference is that in our situation we can not apply
auxiliary results for countable complexes which were used
in \cite{cv}.

Consider a complex $L \in \Lambda^{0}_{4}$. By definition of
$\Lambda_{4}^{0}$, there is a metrizable compactum $C$ such that
\begin{equation}\label{E:3.1}
\operatorname{e-dim}C = \left[ L\right] .
\end{equation}
Set $M = M_{C}^{4}$, where $M_{C}^{4}$ is a manifold
from Theorem \ref{T:3.1}. We claim that this is a required manifold.
First of all, $M$ is countably compact. Hence, in view of Proposition
\ref{P:2.11},
\begin{equation}\label{E:3.2}
\operatorname{e-dim}M = \operatorname{e-dim}\beta M .
\end{equation}
Further, by Corollary \ref{C:2.7} and Theorem \ref{T:3.1}, we have
\begin{equation}\label{E:3.3}
\operatorname{e-dim}\beta M \geq \operatorname{e-dim}
(\beta M -M) \geq \operatorname{e-dim}C = \left[ L\right] .
\end{equation}
Now we apply Proposition \ref{P:2.12} to the pair
$\left( S_{\varphi}^{3},C\right)$. Since $S_{\varphi}^{3}-C$ is
open in $S^{3}$, Theorem \ref{T:2.22} yields
\begin{equation}\label{E:3.4}
\operatorname{e-dim}S_{\varphi}^{3} \leq
\max\left\{ \left[ S^{3}\right] ,\left[ L\right] \right\} =
\left[ L\right] .
\end{equation}
Finally, let us apply Proposition \ref{P:2.12} to the pair
$\left( \beta M, S_{\varphi}^{3}\right)$. Since $\dim Y \leq 4$
for any compactum $Y \subseteq M$, by Theorem \ref{T:2.22} and
inequality (\ref{E:3.4}), we obtain
\begin{equation}\label{E:3.5}
\operatorname{e-dim}\beta M \leq \max \left\{
\left[ S^{n}\right] ,\left[ L\right] \right\} = \left[ L\right] .
\end{equation}
Inequalities (\ref{E:3.3}) and (\ref{E:3.5}) yield
\[ \operatorname{e-dim}\beta M = \left[ L\right] .\]
Thus, equality (\ref{E:3.2}) finishes the proof of
Theorem \ref{T:3.2}.
\end{proof}

As corollaries of Theorem \ref{T:3.2} we discuss several examples of
complexes $L \in \Lambda_{4}^{0}$ with certain curious properties.
First of all we recall two results needed for our discussion.

\begin{pro}\cite[Proposition 2.6]{hand}\label{P:3.3}
Let $K$ be an $n$-dimensional locally compact polyhedron.
Then $\operatorname{e-dim}K = \left[ S^{n}\right]$.
\end{pro}

\begin{pro}\cite[Corollary 2.3]{hand}\label{P:3.4}
Let $L$ be a simplicial complex homotopy dominated by a
finite complex. Then there exists a metrizable compactum
$X^{L}$ such that $\operatorname{e-dim}X^{L} = \left[ L\right]$.
\end{pro}

\begin{rem}\label{R:3.5}
{\em It follows from the proof of Theorem \ref{T:3.2}
that for every complex $L \in \Lambda_{4}^{0}$ there exists
a metrizable compactum $C_{L}$ such that}
\begin{equation}\label{E:3.6}
\operatorname{e-dim}M^{4,L} = \operatorname{e-dim}
M_{C_{L}}^{4} = \operatorname{e-dim}C_{L} = \left[ L\right] .
\end{equation}
\end{rem}

\begin{ex}\label{Ex.3.6}
{\em Let ${\mathcal L} = \{ S^{n} \colon n \geq 4\}$ and $C_{n} = I^{n}$.
Then from (\ref{E:3.6}) and Proposition \ref{P:3.3} we obtain that
$\dim M^{4,S^{n}} = n$ -- the fact proved earlier in \cite{vvf1}.}
\end{ex}

\begin{defin}\label{D:3.7}
{\em Let $L_{n} = M(\text{\f Z}_{2}, n+1) \vee S^{n+1}$,
where $M(\text{\f Z}_{2}, n+1)$ is the} Moore complex,
{\em i.e. the complex obtained from $(n+1)$-dimensional
disk $B^{n+1}$ by attaching to its boundary $S^{n}$
the disk $B^{n+1}$ via the map $S^{n} \to S^{n}$ of degree $2$.
It is clear that $L_{n}$ is a finite complex with
$\left[ S^{n} \right] < \left[ L_{n}\right] <
\left[ S^{n+1}\right]$.}
\end{defin}

\begin{ex}\label{Ex.3.8}
{\em Let ${\mathcal L} = \{ L_{n} \colon n \geq 4\}$ and let
$C_{n}$ be a metrizable compactum with
$\operatorname{e-dim}C_{n} = \left[ L_{n}\right]$
(see Proposition \ref{P:3.4}). Then
$\operatorname{e-dim}M^{4,L_{n}} = \left[ L_{n}\right]$.}
\end{ex}

\begin{cor}\label{C:3.9}
Assuming $\lozenge$, there exists a differentiable,
countably compact,
perfectly normal, hereditarily separable
$4$-manifold $M^{4}$ such that
$\left[ S^{4}\right] < \operatorname{e-dim}M^{4}
< \left[ S^{5}\right] $.
\end{cor}


\section{On the dimension of products of manifolds}
The Stone-\v{C}ech remainder $\beta X - X$ of a space $X$ is
denoted by $X^{\ast}$.
\begin{lem}\label{L:4.1}
Let $M_{i}$ be a countably compact $n_{i}$-manifold of dimension
$\dim M_{i} = m_{i}$, $i = 1,2$. Then
\[ \dim (M_{1}\times M_{2}) = \max\{ n_{1}+m_{2}, n_{2}+m_{1},
\dim(M_{1}^{\ast}\times M_{2}^{\ast})\} .\]
\end{lem}
\begin{proof}
Because each manifold is a $k$-space (being first countable)
it follows from \cite[Theorem 3.10.13]{engelking1} that
$M_{1}\times M_{2}$ is countably compact. Hence, by
Gliksberg's theorem \cite{gli}, $M_{1}\times M_{2}$ is pseudocompact
and $\beta (M_{1}\times M_{2}) = \beta M_{1}\times \beta M_{2}$.
By Proposition \ref{P:2.1},
$\dim (M_{1}\times M_{2}) = \dim \beta (M_{1}\times M_{2})$.
So we have to find out the exact value of
$\dim \beta (M_{1}\times M_{2})$. 

Let $X = \beta M_{1}\times \beta M_{2}$,
$F = M_{1}^{\ast}\times M_{2}^{\ast}$ and $U = X-F$.
By Dowker's theorem
\cite{dowk},
\begin{equation}\label{E:4.1}
\dim X = \max\{ \dim F , k\} ,
\end{equation}
where
\begin{equation}\label{E:4.2}
k = \sup\{ \dim Y \colon Y \subseteq U,\; Y \;
\text{is closed in}\; X\} .
\end{equation}
It is clear, that each $Y$ from (\ref{E:4.2}) is contained
in some $Y^{\prime} = (K_{1}\times\beta M_{2}) \cup
(\beta M_{1}\times K_{2})$, where $K_{i} \subseteq M_{i}$
is a finite sum of $n_{i}$-dimensional cubes, $i = 1,2$.
By Morita's theorem \cite{morita},
\begin{equation}\label{E:4.3}
\dim (K \times Z) = \dim K + \dim Z ,
\end{equation}
whenever $Z$ is a paracompact space and $K$ is a compact
polyhedron. By the finite sum theorem for $\dim$, (\ref{E:4.3})
yields
\begin{equation}\label{E:4.4}
\dim (K_{i} \times \beta M_{j}) = \dim K_{i}+\dim\beta M_{j} .
\end{equation}
Consequently,
\begin{multline*}
\dim Y^{\prime} = \max\{ \dim (K_{1}\times\beta M_{2}),
\dim (\beta M_{1}\times K_{2})\}
\stackrel{\text{by} (\ref{E:4.4})}{=}\\
\max\{ n_{1}+\dim \beta M_{2}, \dim\beta M_{1}+n_{2}\}
\stackrel{\text{by Proposition} \ref{P:2.1}}{=}\\
\max\{ n_{1}+\dim M_{2}, \dim M_{1}+n_{2}\} =
\max\{ n_{1}+m_{2}, m_{1}+n_{2}\} .
\end{multline*}
Equality (\ref{E:4.1}) finishes the proof of Lemma \ref{L:4.1}.
\end{proof}

\begin{cor}\label{C:4.2}
Let $M_{i}$ be a countably compact $n_{i}$-manifold of
dimension $\dim M_{i} = m_{i}$, $i = 1,2$. If
$\max\{ n_{1}+m_{2}, n_{2}+m_{1}\} \leq
\dim (M_{1}^{\ast}\times M_{2}^{\ast})$, then
$\dim (M_{1}\times M_{2}) = \dim (M_{1}^{\ast}\times
M_{2}^{\ast})$.
\end{cor}
\begin{pro}\label{P:4.3}
Let $M_{1}$ and $M_{2}$ be countably compact manifolds. Then
\[ \dim (M_{1} \times M_{2}) \leq \dim M_{1}+\dim M_{2} .\]
\end{pro}
\begin{proof}
According to Lemma \ref{L:4.1}, we only need to check that
$\dim (M_{1}^{\ast}\times M_{2}^{\ast}) \leq m_{1}+m_{2}$.
But for any compact spaces $X_{1}$ and $X_{2}$ we have
(see \cite{hemm})
\[ \dim (X_{1}\times X_{2}) \leq \dim X_{1} +\dim X_{2} .\]
Hence, 
\begin{multline*}
\dim (M_{1}^{\ast}\times M_{2}^{\ast}) \leq
\dim M_{1}^{\ast} +\dim M_{2}^{\ast} \leq
\dim\beta M_{1}+\dim\beta M_{2} =\\
 \dim M_{1}+\dim M_{2} =
m_{1}+m_{2} .
\end{multline*}
Proposition \ref{P:4.3} is proved.
\end{proof}

The next statement is an immediate corollary of the countable
sum theorem for Lebesgue dimension.
\begin{lem}\label{L:4.4}
Let $X_{i}$ be metrizable compacta, $i = 1,2$. Let $F_{i}$
be a closed subset of $X_{i}$, and let $U_{i} = X_{i} - F_{i}$.
Then
\begin{multline*}
 \dim (X_{1}\times X_{2}) =\\
 \max\{ \dim (U_{1}\times U_{2}),
\dim (U_{1}\times F_{2}), \dim (F_{1}\times U_{2}),
\dim (F_{1}\times F_{2})\} .
\end{multline*}
\end{lem}

\begin{thm}\label{T:4.5}
Let $m$ be a natural number such that $m \geq 5$. Then, assuming
$\lozenge$, there exists a differentiable, countably compact,
perfectly normal, hereditarily separable $4$-manifold
$M = M_{m}$ such that $\dim M = m$ and
$\dim (M \times M) = 2m-1 < 2\dim M$.
\end{thm}
\begin{proof}
Let $B$ be a two-dimensional metrizable compactum such that
$\dim (B\times B) = 3$. Such a compactum was
constructed by V.~G.~Boltyanski \cite{bolt}. Let
$C = B \times I^{m-2}$.
Then in accordance with (\ref{E:4.3}),
\begin{equation}\label{E:4.5}
\dim C = m ,
\end{equation}
\begin{equation}\label{E:4.6}
\dim (C\times C) = 2m-1 .
\end{equation}
Let $M = M_{C}^{4}$ be a manifold from Theorem \ref{T:3.1}.
We claim that $M$ is a required manifold. Indeed,
by the properties of $M_{C}^{4}$, the set $M^{\ast}-C = U$ is
homeomorphic to an open subset of $S^{3}$. Consequently,
Lemma \ref{L:4.1} and (\ref{E:4.6}) yield
$\dim (M^{\ast}\times M^{\ast}) = 2m-1$. In this situation
Corollary \ref{C:4.2} finishes the proof of theorem \ref{T:4.5}.
\end{proof}

\begin{q}\label{Q:4.6}
{\em Does there exist a $4$-manifold $M$ such that}
\[2\dim M -\dim (M\times M) \geq 2 \;\;?\]
\end{q}
A similar question about two {\em different} manifolds
has a positive solution.
\begin{thm}\label{T:4.7}
Let $m_{1}, m_{2}$ and $r$ be natural numbers such that
$5 \leq m_{1}\leq m_{2}$ and $4+m_{2} \leq r < m_{1}+m_{2}$.
Then assuming $\lozenge$ there exist differentiable, countably
compact, perfectly normal, hereditarily separable $4$-manifolds
$M_{1}$ and $M_{2}$ of dimension $\dim M_{i} = m_{i}$ such that
\[ \dim (M_{1}\times M_{2}) = r < m_{1}+m_{2} =
\dim M_{1}+\dim M_{2} .\]
\end{thm}
\begin{proof}
We follow the proof of Theorem \ref{T:4.5}. First let us recall
the following result \cite[\S 2, Corollary 2]{dra}.
\begin{itemize}
\item
{\em For all natural numbers $m_{1}, m_{2}$ and $r$ such that
$m_{1} \leq m_{2}$ and $m_{2} < r \leq m_{1}+m_{2}$, there
exist metrizable compacta $X_{1}$ and $X_{2}$ such that
$\dim X_{i} = m_{i}$
and $\dim (X_{1}\times X_{2}) = r$.}
\end{itemize}
Set $M_{i} = M_{X_{i}}^{4}$, where $X_{1}$ and $X_{2}$ are the
above mentioned compacta with $m_{1}, m_{2}$ and $r^{\prime}$
satisfying
inequalities $m_{2} < r^{\prime} \leq m_{1}+m_{2}$.
From Lemma \ref{L:4.4}
we get $\dim (M_{1}^{\ast}\times M_{2}^{\ast}) =
\max\{ 3+m_{2}, r^{\prime}\}$. On the other hand, for $k$
from (\ref{E:4.2}) we have $k = \dim Y^{\prime} = 4 +m_{2}$.
In view of Lemma
\ref{L:4.1} we have $\dim(M_{1}\times M_{2}) =
\max\{ 4 + m_{2}, r^{\prime}\}$. Theorem \ref{T:4.7} is proved.
\end{proof}

\begin{rem}\label{R:4.8}
{\em As we have seen the dimension of the product  of manifolds
can be much less than the sum of their dimensions. But, since
our manifolds $M_{i}$ are countably compact,
by Proposition \ref{P:4.3} we have $\dim (M_{1}\times M_{2})
\leq \dim M_{1}+\dim M_{2}$.}
\end{rem}

\begin{q}\label{Q:4.9}
{\em Are there manifolds $M_{1}$ and $M_{2}$ such that}
\[\dim (M_{1}\times M_{2})
> \dim M_{1}+\dim M_{2}\;\; ?\]
\end{q}

\section{On the dimension of subsets of $4$-manifolds}

\begin{thm}\label{T:5.1}
Assuming $\lozenge$, there exists an infinite-dimensional,
differentiable, countably compact, perfectly normal,
hereditarily separable $4$-manifold $M^{4}$ such that
for every closed set $F \subseteq M^{4}$ we have
\[ \text{either}\;\; \dim F \leq 4\;\;\text{or}\;\;
\dim F = \infty .\]
\end{thm}
\begin{proof}
Let $M^{4} = M_{C}^{4}$, where $M_{C}^{4}$ is a manifold
from Theorem \ref{T:3.1} and $C$ is well known Henderson's
infinite-dimensional compactum with no
positive-dimensional compact subsets \cite{hend}. By
Proposition \ref{P:2.1} and Theorem \ref{T:3.1},
$\dim M^{4} = \infty$. Now let $F$ be a closed subset
of $M^{4}$ such that $\dim F \geq 5$. Then in
view of Proposition \ref{P:2.1},
\[ \dim \beta F \geq 5 .\]
But, since $M^{4}$ is normal, $\beta F =
\operatorname{cl}_{\beta M^{4}}F$. Set 
\[ A =
\left( \operatorname{cl}_{\beta M^{4}}F\right) \cap
\left( \beta M^{4} - M^{4}\right) .\]
By Proposition
\ref{P:2.12} for the pair $(\beta F, A)$, we have
$\dim A \geq 5$. But in accordance with Theorem
\ref{T:3.1},
$\beta M^{4} - M^{4}$ is a disjoint sum of $C$ and some
open subset of $S^{3}$. Hence,
\[ \dim (A \cap C) \geq 5 .\]
Therefore, by the property of Henderson's compactum,
$\dim (A \cap C) = \infty$. Consequently, 
\[ \dim F = \dim \beta F \geq \dim C \geq \dim
(A \cap C) = \infty .\]
Theorem \ref{T:5.1} is proved.
\end{proof}

The next statement is a generalization of Theorem
\ref{T:3.2}.
\begin{thm}\label{T:5.2}
Let ${\mathcal L}$ be a countable subset of
$\Lambda_{4}^{0}$ (see Theorem \ref{T:3.2}). Then,
assuming $\lozenge$, there exists a differentiable,
countably compact, perfectly normal, hereditarily
separable $4$-manifold $M^{4}$ which admits a family
$U_{L}$, $L \in {\mathcal L}$, of open subsets  of
extension dimension
$\operatorname{e-dim}U_{L} = \left[ L\right]$.
Moreover, one can choose sets $U_{L}$ in such a way that
\begin{itemize}
\item[(1)]
either $\operatorname{cl}(U_{L}) \cap
\operatorname{cl}(U_{L^{\prime}}) = \emptyset$ if
$L \neq L^{\prime}$,
\item[(2)]
or $\cap\{ U_{L} \colon L \in {\mathcal L}\} \neq \emptyset$.
\end{itemize}
\end{thm}
\begin{proof}
By Theorem \ref{T:3.2}, there is a metrizable compactum
$X_{L}$ of extension dimension
$\operatorname{e-dim}X_{L} = \left[ L\right]$. Let $C$
be the Alexandroff compactification of the discrete sum
$\oplus\{ X_{L} \colon l \in {\mathcal L}\}$ of these compacta.
Now we define $M^{4}$ as a manifold $M_{C}^{4}$ from
Theorem \ref{T:3.1}. Since $\{ X_{L} \colon L \in {\mathcal L}\}$
is a discrete family in a compact space
$\beta M^{4}$, there is a disjoint family of
neighborhoods $V_{L}$ of $X_{L}$ in $\beta M^{4}$.
We may assume also, that
\[ \operatorname{cl}_{\beta M^{4}}(V_{L}) \cap
\operatorname{cl}_{\beta M^{4}}(V_{L^{\prime}}) =
\emptyset\;\text{if}\; L \neq L^{\prime} .\]
Now, in case (1), we set
\begin{equation}\label{E:5.1}
U_{L} = V_{L} \cap M^{4} .
\end{equation}
To realize case (2) we take an open metrizable subset
$U \subseteq M^{4}$ and set
\begin{equation}\label{E:5.2}
U_{L} = (U \cup V_{L})\cap M^{4} .
\end{equation}

{\em Claim}. $\operatorname{e-dim}
(\operatorname{cl}_{M^{4}}U_{L}) = \left[ L\right]$.

{\em Proof}. By Proposition \ref{P:2.11}, it suffices
to verify that
\begin{equation}\label{E:5.3}
\operatorname{e-dim}\beta \left
(\operatorname{cl}_{M^{4}}U_{L}\right) = \left[ L\right] .
\end{equation}
But
\begin{equation}\label{E:5.4}
\beta\left( \operatorname{cl}_{M^{4}}U_{L} \right)
= \operatorname{cl}_{\beta M^{4}}\left
(\operatorname{cl}_{M^{4}}U_{L}\right) =
\operatorname{cl}_{\beta M^{4}}\left( U_{L}\right) .
\end{equation}
On the other hand, according to (\ref{E:5.1}), $U_{L}$ is
dense in $V_{L}$. Hence,
\begin{equation}\label{E:5.5}
\operatorname{cl}_{\beta M^{4}}\left( U_{L}\right) =
\operatorname{cl}_{\beta M^{4}}\left( V_{L}\right) .
\end{equation}
Let $\Phi_{L} = M^{4} \cap
\operatorname{cl}_{\beta M^{4}}\left( V_{L}\right)$
and let $F_{L} = \operatorname{cl}_
{\beta M^{4}}\left( V_{L}\right) - M^{4}$. Since
$\Phi_{L}$ is dense in
$\operatorname{cl}_{\beta M^{4}}\left(V_{L}\right)$, we have
\begin{equation}\label{E:5.6}
\operatorname{cl}_{\beta M^{4}}\Phi_{L} =
\Phi_{L} \cup F_{L} = \operatorname{cl}_
{\beta M^{4}}\left( V_{L}\right) .
\end{equation}
Hence,
\begin{equation}\label{E:5.7}
\beta\Phi_{L} = \Phi_{L} \cup F_{L} .
\end{equation}
For every compactum $Y \subseteq \Phi_{L}$ we have
$\dim Y \leq 4$. On the other hand, $F_{L}
\supseteq X_{L}$ and $F_{L} \cap X_{L^{\prime}} =
\emptyset$ whenever $L \neq L^{\prime}$. Hence,
$\operatorname{e-dim}F_{L} = \left[ L\right]$. Now
we apply proposition \ref{P:2.12} to the pair
$\left( \beta \Phi_{L}, F_{L}\right)$. We have
\begin{equation}\label{E:5.8}
\operatorname{e-dim}\beta\Phi_{L} =
\operatorname{e-dim}F_{L} = \left[ L\right] .
\end{equation}
Finally, conditions (\ref{E:5.7}), (\ref{E:5.6}),
(\ref{E:5.5}) and (\ref{E:5.4}) give us the required
equality (\ref{E:5.3}). This finishes proof of our Claim.

In order to prove the equality
$\operatorname{e-dim}U_{L} = \left[ L\right]$ we
need more general
version of Theorem \ref{T:2.21}. Its proof is based on
the fact that a countable simplicial complex is an
$ANE$ for the class of all normal spaces.

\begin{itemize}
\item
Suppose that a normal space $X$ is the union
of its closed countably compact subsets $X_{i}$,
$i \in \omega$. If $\operatorname{e-dim}X_{i} \leq
\left[ L\right]$ for each $i \in \omega$, then
$\operatorname{e-dim}X \leq \left[ L\right]$.
\end{itemize}

In order to finish the proof of Theorem \ref{T:5.2} represent
$U_{L}$ as the union of an increasing sequence
\[ U_{L}^{0} \subseteq \operatorname{cl}_{M^{4}}U_{L}^{0}
\subseteq U_{L}^{1} \subseteq \cdots \]
and apply our claim. Theorem \ref{T:5.2} is proved.
\end{proof}


\begin{thebibliography}{99}



\bibitem{alpas}
P.~S.~Alexandrov and B.~A.~Pasynkov, {\em
Introduction to the Dimension Theory}, ``Nauka",
Moscow, 1973.

\bibitem{bolt}
V.~G.~Boltyanski, {\em An example of a two-dimensional
compactum whose topological square is three-dimensional},
Dokl. Akad. Nauk SSSR {\bf 67} (1949), 597--599.

\bibitem{coho}
A.~Chigogidze, {\em Cohomological dimension of Tychonov spaces},
Topology Appl. {\bf 79} (1997), 197--228.

\bibitem{hand}
\bysame, {\em Infinite dimensional topology and shape theory},
Handbook of Geometric Topology (edited by R.Daverman and R.Sher), to appear.

\bibitem{cv}
A.~Chigogidze and V.V.Fedorchuk, {\em On extension
dimension of non-metrizable manifolds}, Vestnik Moscow Univ., submitted.


\bibitem{dowk}
C.~H.~Dowker, {\em Local dimension of normal spaces}, Quart. J.
Math. {\bf 6}:22 (1955), 101--120.

\bibitem{dra}
A.~N.~Dranishnikov, {\em Homological dimension theory},
Uspekhi Mat. Nauk {\bf 43}:4 (1988), 11--55.

\bibitem{dra1}
\bysame, {\em The Eilenberg-Borsuk theorem for maps in an
arbitrary complex}, Mat. Sbornik {\bf 185} (1994), 81--90.

\bibitem{dra2}
\bysame, {\em Extension of mappings into CW-complexes},
Mat. Sbornik {\bf 182} (1991), 47--56.

\bibitem{dra3}
\bysame, {\em On extension theory of compact spaces},
Topology Appl. (the same issue).

\bibitem{engelking1}
R.~Engelking, {\em General Topology}, PWN, Warszaw, 1977.


\bibitem{vvf1}
V.~V.~Fedorchuk, {\em A differentiable manifold with non-coinciding
dimensions},
Topology Appl. {\bf 54} (1993), 221--239.


\bibitem{vvf}
\bysame, {\em On the transfinite and cohomological
dimensions of $4$-manifolds}, Proc. Steklov Math. Inst.
{\bf 211} (1996), 184--202.


\bibitem{fri}
R.~Fritch and R.~A.~Piccinini, {\em Cellular
Structures in Topology},
Cambridge Univ. Press, Cambridge, 1990.

\bibitem{gli}
I.~Gliksberg, {\em Stone-\v{C}ech compactification
of products},
Trans. Amer. Math. Soc. {\bf 90} (1959), 369--382.

\bibitem{hemm}
E.~Hemmingsen, {\em Some theorems in dimension theory for
normal Hausdorff spaces}, Duke Math. J. {\bf 13}
(1946), 495--504.

\bibitem{hend}
D.~W.~Henderson, {\em An infinite-dimensional compactum with no
positive-dimensional compact subsets -- a simpler construction},
Amer. J. Math. {\bf 89} (1967), 105--121.


\bibitem{morita}
K.~Morita, {\em On the dimension of product spaces}, Amer. J.
Math. {\bf 75} (1953), 205--223.

\bibitem{pears}
A.~R.~Pears, {\em Dimension Theory of General Spaces},
Cambridge Univ. Press, Cambridge, 1975.
\end{thebibliography}
\end{document}